# ON THE BENJAMINI–HOCHBERG METHOD

By J. A. Ferreira[1] and A. H. Zwinderman

*University of Amsterdam*

We investigate the properties of the Benjamini–Hochberg method for multiple testing and of a variant of Storey's generalization of it, extending and complementing the asymptotic and exact results available in the literature. Results are obtained under two different sets of assumptions and include asymptotic and exact expressions and bounds for the proportion of rejections, the proportion of incorrect rejections out of all rejections and two other proportions used to quantify the efficacy of the method.

**1. Introduction.** Let $X = \{X_1, X_2, \ldots, X_m\}$ be a set of $m$ random variables defined on a probability space $(\Omega, \mathcal{F}, P)$ such that, for some positive integer $m_0 \leq m$, each of $X_1, X_2, \ldots, X_{m_0}$ has distribution function (d.f.) $F$ and $X_{m_0+1}, \ldots, X_m$ all have d.f.'s different from $F$, and consider the problem of choosing a set $\mathcal{R} \subseteq X$ in such a way that the random variable (r.v.)

$$\Pi_{1,m} = \frac{S_m}{R_m \vee 1},$$

where $R_m = \#\mathcal{R}$ and $S_m = \#(\mathcal{R} \cap \{X_1, \ldots, X_{m_0}\})$, is guaranteed to be small in some probabilistic sense. In more ordinary language, the problem is that of discovering observations in $X$ which do not have d.f. $F$ without incurring a high *proportion of incorrect rejections*—the proportion $\Pi_{1,m}$ of rejected observations which in fact come from $F$.

Benjamini and Hochberg [2] have proposed a method of choosing $\mathcal{R}$ specifically aimed at discovering r.v.'s taking values in the interval $[0,1]$ that tend to be smaller than standard uniform r.v.'s and which, given $\delta > 0$, guarantees that $E(\Pi_{1,m}) \leq \delta$ under certain conditions. The method consists of

Received February 2005; revised August 2005.

[1]Supported in part by "Nederlandse Hartstichting," the Dutch Heart Foundation, through Grant 2003B215.

fixing $q \in [0, 1]$, computing

$$(1.1) \qquad R_m = \max\left\{i : X_{i:m} \leq q\frac{i}{m}\right\},$$

where $0 \leq X_{1:m} \leq \cdots \leq X_{m:m} \leq 1$ denote the order statistics of $X$, and setting $\mathcal{R} = \{X_{1:m}, \ldots, X_{R_m:m}\}$. In its simplest form, the Benjamini–Hochberg theorem states that if $\mathcal{R}$ is chosen according to this procedure and $X_1, X_2, \ldots, X_{m_0}$ are independent and standard uniform and independent of $X_{m_0+1}, \ldots, X_m$, then $E(\Pi_{1,m}) = q\gamma$, where $\gamma := m_0/m$, a property usually expressed by saying that the Benjamini–Hochberg procedure controls the *false discovery rate*—the number $E(\Pi_{1,m})$.

The Benjamini–Hochberg procedure seems somewhat mysterious from (1.1) alone; an explanation as to why it does work in the appropriate circumstances will be given below.

Benjamini and Hochberg [2] formulated their ideas in the context of *multiple testing*. Here, rejecting observations in $X$ is interpreted as rejecting hypotheses among $m$ null hypotheses $H_0^1, \ldots, H_0^m$, of which only the first $m_0$ are true, on the basis of *p-values* $X_1, \ldots, X_m$ that result from the observation of certain test statistics $Y_1, \ldots, Y_m$. Although the hypotheses tested may be arbitrary, the *p*-values are assumed to be given by $X_i = 1 - F_i(Y_i)$, where $F_i$ is the d.f. of $Y_i$ under $H_0^i$; furthermore, in the most general case considered by Sarkar [15] $X_1, X_2, \ldots, X_{m_0}$ need not be independent and are only assumed to be *sub-uniform* in the sense that $P(X_i \leq x) \leq x$ for all $x \in [0, 1]$. [*Note*: In general, $P(X_i \leq x) \geq x$, rather than $P(X_i \leq x) \leq x$: If $F$ is a d.f. and $F^{-1}(u) = \min\{t : F(t) \geq u\}$ then $F(t) \geq u \Leftrightarrow t \geq F^{-1}(u)$, and $F(F^{-1}(u)-) \leq u$; therefore, $P(X_i \leq x) = P(F_i(Y_i) \geq 1 - x) = P(Y_i \geq F_i^{-1}(1-x)) = 1 - F_i(F_i^{-1}(1-x)-) \geq x$ with equality for all $x$ if and only if $F_i$ is continuous. Thus (see, e.g., the proof of Theorem 2.1), under the assumptions usually made in the literature, the Benjamini–Hochberg theorem actually states that $E(\Pi_{1,m}) \geq q\gamma$. If the method is modified by using strict inequality in (1.1) and the *p*-values are defined by $X_i = F_i(Y_i)$ (which represents no loss of generality), then $E(\Pi_{1,m}) \leq q\gamma$ with equality if $Y_1, \ldots, Y_{m_0}$ are continuous, because $P(X_i < x) = P(F_i(Y_i) < x) = P(Y_i < F_i^{-1}(x)) = F_i(F_i^{-1}(x)-) \leq x$.]

Most common multiple testing procedures tend to be either too conservative or too liberal—they either miss the chance of detecting many false hypotheses in the fear of incorrectly rejecting one hypothesis (the case of the Bonferroni method), or they incur a very large proportion of false positives in the greed of finding significant results (the case of "uncritical testing," in which all hypotheses yielding *p*-values below $q$, say, are rejected). Benjamini and Hochberg's [2] motivation in proposing to control the false discovery rate was to achieve a balance between these two extremes: in many problems—especially in those involving many hypotheses—it is acceptable to incorrectly reject some hypotheses as long as they make up only a small proportion of all the hypotheses rejected; and allowing for this proportion of false positives yields a substantial proportion of true discoveries. We were led to the Benjamini–Hochberg approach to multiple testing by considering one such problem: "gene discovery" in the context of heart disease, where the objective is to discover genetic variables which determine or influence a number of phenotypical variables. "Gene expression" studies provide other examples of problems where the control of the false discovery rate is important; see, for example, Tusher, Tibshirani and Chu [22], Dudoit, Schaffer and Boldrick [7], Reiner, Yekutieli and Benjamini [14], Fan et al. [8] and McLachlan, Do and Ambroise [12]. Some of these authors actually use variants of the Benjamini–Hochberg method based on estimating the proportion of incorrect rejections out of all rejections that result from rejecting all $p$-values below $t$ as a function of $t$, a procedure which for $t = qR_m/m$ is equivalent to Benjamini and Hochberg's.

As outlined in our first paragraph, the problem of choosing $\mathcal{R}$ in a way that controls $\Pi_{1,m}$ seems to arise in other contexts as well. For instance, in data analyses of "contaminated" data, where a majority of elements form a sample from some population but a minority do not, $\mathcal{R}$ records those observations thought to be "outliers," and it is naturally of interest to seek a choice of $\mathcal{R}$ that keeps $\Pi_{1,m}$ small so that not too many of the good observations are thrown away. In the more general formulation, the variables $X_{m_0+1}, \ldots, X_m$ need not behave in a more extreme way than $X_1, \ldots, X_{m_0}$; they simply have d.f.'s that differ from $F$, and the problem, then, can be further translated into that of identifying a mixture of two populations given the knowledge of the law describing one of them. This is a useful point of view in that it helps us to put the Benjamini–Hochberg method into a context of goodness of fit, which is not just more general but also illuminating as far as the workings and the limitations of the method are concerned. More specifically, the problem could, in principle, be solved by choosing $\mathcal{R}$ as the subset of $X$ for which a goodness of fit test of $F$ performed with $X \setminus \mathcal{R}$ yields the smallest discrepancy among the discrepancies based on all subsets of $X$. As we shall see, what the Benjamini–Hochberg method does is just this, except that the subsets considered are of the form $\{X_{1:m}, \ldots, X_{r:m}\}$ for some $r$.

Let $H_m$ denote the empirical d.f. of $X$; then (the second identity here is known and has been used before in this context; e.g., see [1] and [9])

$$\{R_m \geq r\} = \bigcup_{k=r}^{m} \left\{ X_{k:m} \leq q\frac{k}{m} \right\}$$

$$= \bigcup_{k=r}^{m} \left\{ \sum_{j=1}^{m} \mathbb{1}_{\{X_j \leq qk/m\}} \geq k \right\}$$

$$(1.2) \quad \begin{aligned} &= \bigcup_{k=r}^{m} \left\{ H_m\left(q\frac{k}{m}\right) \geq \frac{k}{m} \right\} \\ &= \bigcup_{k=r}^{m} \left\{ \frac{H_m(qk/m) - qk/m}{qk/m} \geq \frac{1-q}{q} \right\} \\ &= \left\{ \max_{t=qr/m,\ldots,q(m-1)/m,q} \frac{H_m(t) - t}{t} \geq \frac{1-q}{q} \right\}, \end{aligned}$$

$r = 0, 1, \ldots, m$, so the procedure rejects the $r$ lower order statistics if and only if

$$\max_{t=qr/m,\ldots,q(m-1)/m,q} \frac{H_m(t) - t}{t} \geq \frac{1-q}{q}$$

and

$$\max_{t=q(r+1)/m,\ldots,q(m-1)/m,q} \frac{H_m(t) - t}{t} < \frac{1-q}{q}.$$

In other words, the $r$ lower order statistics are rejected whenever the goodness of fit statistics

$$(1.3) \quad \max_{t=qk/m,\ldots,q(m-1)/m,q} \frac{H_m(t) - t}{t} \approx \max_{t \in [qk/m,q]} \frac{H_m(t) - t}{t},$$

$k = 1, \ldots, m$, indicate a relatively big discrepancy between $H_m$ and the uniform d.f. over $[qr/m, q]$, and a relatively small one over $[q(r+1)/m, q]$, indicating that most of the nonuniform observations lie in the interval $(0, qr/m]$; the standard for comparison, $(1-q)/q$, corresponds to the biggest discrepancy of $(H_m(t) - t)/t$ one could get at $t = q$, and the choice of $q$ determines the interval $(0, q]$ to be "scanned" for discrepancies.

The function on the right-hand side in (1.3) is Rényi's statistic, a well-known goodness of fit statistic for testing the uniform distribution; it is a one-sided statistic of the Kolmogorov–Smirnov type, devised to detect distributions with too much mass in the lower tail, scaled by the standard uniform distribution in order to inflate the discrepancies that occur at lower values.

From the version of the "ballot theorem" given on page 113 of [11], we know that if $X_1, \ldots, X_m$ are independent standard uniform r.v.'s, then $P(H_m(t) \leq t/q \ \forall t \in (0, q]) = 1 - q$ for all $m \in \mathbb{N}$ and $q \in [0, 1]$, from which it follows that the probability that the Benjamini–Hochberg method yields no rejections satisfies $P(R_m < 1) \sim 1 - P(\sup_{0 < t \leq q}(H_m(t) - t)/t \geq (1-q)/q) = 1 - q$. Thus, if the hypothesis that the variables are a standard uniform random sample is taken as the null and the type I error is defined as the incorrect rejection of at least one $p$-value, $q$ can be interpreted as the approximate significance level. (We thank a referee for posing a question which led to this observation.)

The connection between the Benjamini–Hochberg procedure and goodness of fit has been hinted at by other authors (e.g., [5, 6, 13]), but this seems to be the first explicit link to be exhibited. In their seminal work Benjamini and Hochberg [2] provided some justification of the appropriateness of their method, and so did Storey [18] in connection with one of the variants mentioned above; the present explanation provides further insight into the workings of the method, as well as to its domain of applicability.

The objective of this article is to investigate the main properties of the Benjamini–Hochberg method, extending and complementing the results of Benjamini and Hochberg [2], Genovese and Wasserman [9] and Storey, Taylor and Siegmund [19], focusing particularly on its asymptotic aspects as $m \to \infty$, $m_1 := m - m_0 \to \infty$ and $\gamma$ remains fixed. In Section 2 we extend the Benjamini–Hochberg theorem and prove some results on the convergence in probability of $R_m$ to infinity, and of $\Pi_{1,m}$ to $q\gamma$, in what is essentially the setting originally adopted by Benjamini and Hochberg [2]: $X_1, \ldots, X_{m_0}$ are independent and sub-uniform, and independent of $X_{m_0+1}, \ldots, X_m$, but the latter can be anything. This set of assumptions is very asymmetric in that too much is assumed from one set and nothing is assumed from the other, but the results are potentially useful in a number of practical situations. In fact, *the proofs of Section 2 go through if the assumptions just stated hold conditionally on a sigma field $\mathcal{G} \subset \mathcal{F}$*, hence if $X_1, \ldots, X_{m_0}$ are, for each $m_0$, part of an infinite exchangeable sequence independent of $X_{m_0+1}, \ldots, X_m$, and so the results are more general than stated. (See [4] and [15] for the Benjamini–Hochberg theorem under general dependence conditions. Recent parallel developments in this area can be found in [10] and [17].)

But more interesting, perhaps, is that the results proved in Section 2 actually hold in an asymptotic way under the rather general assumptions introduced by Storey, Taylor and Siegmund [19]. These assumptions, which essentially amount to the convergence of the sequence of empirical distributions, are more balanced and seem more realistic. In our work in Sections 3 and 4 we adopt essentially the assumptions of [19] and obtain results which are parallel to theirs, namely about the convergence in probability of $R_m/m$ and $\Pi_{1,m}$; our approach allows some extensions and, we think, the quickest and most transparent treatment of the main properties of the Benjamini–Hochberg method. The results of Section 3 are extended in Section 4 to a slight modification of Storey's [18] generalization of the Benjamini–Hochberg method, whose practical relevance and range of applicability are illustrated by the statements of Theorem 4.1.

Before proceeding, let us introduce two statistical measures often used to assess the performance of the Benjamini–Hochberg method,

$$\Pi_{2,m} = \frac{R_m - S_m}{m - m_0} \equiv \frac{R_m - S_m}{m_1} \quad \text{and} \quad \Pi_{3,m} = 1 - \frac{m_0 - S_m}{(m - R_m) \vee 1}.$$

The first is the *proportion of correctly rejected observations* out of $\{X_{m_0+1}, \ldots, X_m\}$, and its expected value will be called *average power*, or simply *power*; it is the most popular and perhaps most straightforward efficacy measure considered in the literature. The second is the *proportion of incorrect nonrejections among nonrejections* and has been introduced by Genovese and Wasserman [9] as a dual quantity to $\Pi_{1,m}$; its expected value is called *false nondiscovery rate*. The latter seems to be a particularly useful concept in the context of "outlier detection" mentioned above, where one would like to keep only a small number of outliers out of all the observations judged to have come from $F$; in the multiple testing context it seems more difficult to interpret than average power; but see Proposition 2.3 for an interpretation in terms of the Benjamini–Hochberg method.

## 2. Results in the original setting.

Unless stated otherwise, $X_1, \ldots, X_{m_0}$ will be assumed independent and such that $P(X_i \leq x) \leq F(x) := x$ for $x \in [0, 1]$, and independent of $\{X_{m_0+1}, \ldots, X_m\}$. In the sequel, by $X_{i:m-j}^{(j)}$ we shall mean the $i$th order statistic of the set $X^{(j)} := X \setminus \{X_1, \ldots, X_j\}$, $j = 1, \ldots, m_0$, and by $R_m^{(j)}(X^{(j)})$ the number of rejections that result from applying to $X^{(j)}$ the modified form of the Benjamini–Hochberg procedure obtained by replacing $i$ on the right-hand side of the inequalities in (1.1) by $i+j$; we shall also write $R_m = R_m^{(0)}(X)$, $X_{i:m} = X_{i:m}^{(0)}$, $X = X^{(0)}$. By the *standard uniform case*, we mean the case where $X_1, \ldots, X_{m_0}$ are standard uniform r.v.'s.

Our first result gives upper bounds on the moments of $\Pi_{1,m}$ and $S_m$, and contains Benjamini and Hochberg's [2] theorem as a special case.

THEOREM 2.1. *We have*

$$(2.1) \quad E[(\Pi_{1,m})^k] \leq \sum_{j=1}^{k} \left(q \frac{m_0}{m}\right) \cdots \left(q \frac{m_0 - j + 1}{m}\right) E[(j + R_m^{(j)}(X^{(j)}))^{j-k}]$$

*and*

$$(2.2) \quad E(S_m{}^k) \leq \sum_{j=1}^{k} \left(q \frac{m_0}{m}\right) \cdots \left(q \frac{m_0 - j + 1}{m}\right) E[(j + R_m^{(j)}(X^{(j)}))^{j}]$$

*for $k = 1, 2, \ldots, m_0$, the inequalities being achieved for all $q$ only in the standard uniform case.*

PROOF. We only prove (2.1); the proof of (2.2) is very similar. It will be evident that there is no loss of generality in assuming that $X_1, \ldots, X_{m_0}$ have the same distribution. Observe first that, for $0 \leq r \leq m$ (setting $X_{0:m} = 0$), $R_m = r \Leftrightarrow X_{r:m} \leq q\frac{r}{m} \wedge X_{s:m} > q\frac{s}{m} \ \forall s > r$, and that, for $1 \leq r \leq m$,

$$\left\{ X_1 \leq q\frac{r}{m}, R_m = r \right\}$$

$$= \left\{ X_1 \leq q\frac{r}{m}, X_{r:m} \leq q\frac{r}{m}, X_{s:m} > q\frac{s}{m} \ \forall s > r \right\}$$

$$= \left\{ X_1 \leq q\frac{r}{m}, X_{r-1:m-1}^{(1)} \leq q\frac{r}{m}, X_{s-1:m-1}^{(1)} > q\frac{s}{m} \ \forall s > r \right\}$$

$$= \left\{ X_1 \leq q\frac{r}{m}, X_{r-1:m-1}^{(1)} \leq q\frac{r}{m}, X_{s:m-1}^{(1)} > q\frac{s+1}{m} \ \forall s > r-1 \right\}$$

$$= \left\{ X_1 \leq q\frac{r}{m}, R_m^{(1)}(X \setminus \{X_1\}) = r-1 \right\}.$$

Similarly,

$$\left\{ X_1 \leq q\frac{r}{m}, \dots, X_j \leq q\frac{r}{m}, R_m = r \right\}$$

$$= \left\{ X_1 \leq q\frac{r}{m}, \dots, X_j \leq q\frac{r}{m}, R_m^{(j)}(X^{(j)}) = r-j \right\}$$

for $r = j, j+1, \dots, m$, $j = 0, 1, \dots, m_0$. Thus, since $\{X_1, \dots, X_j\}$ and $X^{(j)}$ are independent if $j \leq m_0$, we have

$$E[(\Pi_{1,m})^k] = \sum_{r=1}^{m} E\left( \frac{S_m^k}{r^k} \mathbb{1}_{\{R_m(X)=r\}} \right)$$

$$= \sum_{r=1}^{m} \frac{1}{r^k} E\left[ \left( \sum_{s=1}^{m_0} \mathbb{1}_{\{X_s \leq qr/m\}} \right)^k \mathbb{1}_{\{R_m(X)=r\}} \right]$$

$$= \sum_{r=1}^{m} \sum_{j=1}^{k} \frac{m_0 \cdots (m_0 - j + 1)}{r^k} E[\mathbb{1}_{\{X_1 \leq qr/m, \dots, X_j \leq qr/m\}} \mathbb{1}_{\{R_m(X)=r\}}]$$

$$= \sum_{j=1}^{k} \sum_{r=j}^{m} \frac{m_0 \cdots (m_0 - j + 1)}{r^k}$$
$$\times E[\mathbb{1}_{\{X_1 \leq qr/m, \dots, X_j \leq qr/m, R_m^{(j)}(X^{(j)})=r-j\}}]$$

$$= \sum_{j=1}^{k} \sum_{r=j}^{m} \frac{m_0 \cdots (m_0 - j + 1)}{r^k}$$
$$\times E[\mathbb{1}_{\{X_1 \leq qr/m, \dots, X_j \leq qr/m\}}] E[\mathbb{1}_{\{R_m^{(j)}(X^{(j)})=r-j\}}]$$

$$\leq \sum_{j=1}^{k} \sum_{r=j}^{m} \frac{m_0 \cdots (m_0 - j + 1)}{m^j r^{k-j}} q^j E[\mathbb{1}_{\{R_m^{(j)}(X^{(j)}) = r - j\}}]$$

$$= \sum_{j=1}^{k} \left( q \frac{m_0}{m} \right) \cdots \left( q \frac{m_0 - j + 1}{m} \right) E\left[ \sum_{r=j}^{m} r^{(j-k)} \mathbb{1}_{\{R_m^{(j)}(X^{(j)}) = r - j\}} \right]$$

$$= \sum_{j=1}^{k} \left( q \frac{m_0}{m} \right) \cdots \left( q \frac{m_0 - j + 1}{m} \right) E[(j + R_m^{(j)}(X^{(j)}))^{j-k}],$$

equality holding for all $q$ if and only if $F$ is standard uniform. $\square$

Setting $k = 1$ at each step of the argument yields what is perhaps the simplest and most elementary available proof of the Benjamini–Hochberg theorem; Sarkar [15] gives a proof using similar ideas in a more general setting, and Storey, Taylor and Siegmund [19] give another simple proof based on the optional stopping theorem.

As the following proposition shows, Theorem 2.1 with $k \geq 2$ can be used to derive conclusions about the asymptotic properties of $\Pi_{1,m}$; the proof is given in the Appendix.

PROPOSITION 2.2.  *If $R_m \xrightarrow{P} \infty$, then*

$$(2.3) \qquad \limsup_{m \to \infty} E[(\Pi_{1,m})^k] \leq (q\gamma)^k, \qquad k \in \mathbb{N};$$

*moreover, in the standard uniform case we have*

$$(2.4) \qquad R_m \xrightarrow{P} \infty \quad \textit{if and only if} \quad \Pi_{1,m} \xrightarrow{P} q\gamma.$$

REMARKS.  (i) One practical rule that follows from (2.4) is this: If with large $m$ one rejects a *substantial* (0.1, say, as opposed to 0.001) proportion $R_m/m$ of the sample (indicating $R_m \to \infty$), then one can be sure that $\Pi_{1,m}$, the proportion of incorrect rejections out of all rejections, is not only near, but is practically equal to, the false discovery rate $E(\Pi_{1,m}) = q\gamma$.

(ii) Besides the false discovery rate, some authors consider $E(S_m)/E(R_m \vee 1)$, sometimes called "marginal false discovery rate" (e.g., [20]). When $k = 1$, (2.2) yields $E(S_m)/E[(1 + R_m^{(1)}(X^{(1)})] \leq q\gamma$ with equality in the standard uniform case, which *almost* represents the control of $E(S_m)/E(R_m \vee 1)$. Since, as shown in the proof of Proposition 2.2, $R_m^{(1)}(X^{(1)})$ is asymptotically no smaller than $R_m$, it follows that in the standard uniform case

$$\lim_{m \to \infty} E\left( \frac{S_m}{R_m \vee 1} \right) \equiv q\gamma = \lim_{m \to \infty} \frac{E(S_m)}{1 + E(R_m^{(1)}(X^{(1)}))}$$

$$\leq \liminf_{m \to \infty} \frac{E(S_m)}{1 + E(R_m)} \leq \liminf_{m \to \infty} \frac{E(S_m)}{E(R_m \vee 1)}$$

(an analogous statement with higher moments is also possible).

Because average power is an absolute quantity, there is nothing one can say about it without some information on $X_{m_0+1}, \ldots, X_m$. More precisely, all that one can conclude from Proposition 2.2 is that, because $R_m/m$ can be anything from 0 to 1 (as can be seen from the results of Section 3),

$$\frac{R_m - S_m}{m_1} = \frac{1}{1-\gamma} \frac{R_m}{m} \left(1 - \frac{S_m}{R_m \vee 1}\right)$$

(hence its expected value) is somewhere between 0 and $\frac{1-q\gamma}{1-\gamma} \geq 1$, which, besides the truism that average power is between 0 and 1, only tells us that $R_m/m$ is asymptotically bounded above by $\frac{1-\gamma}{1-q\gamma} \leq 1$.

In contrast, $E(\Pi_{3,m})$, the false nondiscovery rate of Genovese and Wasserman [9], provides a relative measure of the performance of the Benjamini–Hochberg method—it assesses the efficacy of the method in terms of the number of rejections—for which reason one can use a statement like (2.3) to obtain a meaningful upper bound on $\Pi_{3,m}$ (or on its moments):

PROPOSITION 2.3.   *Suppose $\gamma \in (0,1]$. Then*

$$(2.5) \qquad E[(\Pi_{3,m})^l] \leq (1-\gamma)^l + \frac{E[(\Pi_{1,m})^k]}{\gamma^k}, \qquad k, l \in \mathbb{N};$$

*moreover, if $0 \leq q < 1$,*

$$(2.6) \qquad R_m \xrightarrow{P} \infty \quad \Longrightarrow \quad \limsup_{m \to \infty} E[(\Pi_{3,m})^l] \leq (1-\gamma)^l, \qquad l \in \mathbb{N}.$$

PROOF.   If $R_m = 0$, then $\Pi_{3,m} = 1 - \gamma$; if $R_m = m$, then $\Pi_{3,m} = 0$; and if $R_m > 0$, we have $\Pi_{3,m} = 1 - \frac{m\gamma - S_m}{m - R_m} \leq 1 - \gamma \Leftrightarrow \frac{S_m}{R_m} \leq \gamma$. Thus,

$$E[(\Pi_{3,m})^l] = (1-\gamma)^l P(R_m = 0) + E(\Pi_{3,m} \mathbb{1}_{\{S_m/R_m \leq \gamma\}} \mathbb{1}_{\{R_m > 0\}})$$

$$+ E(\Pi_{3,m} \mathbb{1}_{\{S_m/R_m > \gamma\}} \mathbb{1}_{\{R_m > 0\}})$$

$$\leq (1-\gamma)^l P(R_m = 0) + (1-\gamma)^l P(R_m > 0) + E(\mathbb{1}_{\{\Pi_{1,m} > \gamma\}})$$

$$= (1-\gamma)^l + P(\Pi_{1,m} > \gamma) \leq (1-\gamma)^l + \frac{E[(\Pi_{1,m})^k]}{\gamma^k}.$$

By (2.5) and (2.3), $\limsup_{m \to \infty} E[(\Pi_{3,m})^l] \leq (1-\gamma)^l + q^k$, and since $k \in \mathbb{N}$ is arbitrary (2.6) follows.   $\square$

In words, (2.6) says that if $R_m \xrightarrow{P} \infty$, then, asymptotically, the expected proportion of incorrect nonrejections in the Benjamini–Hochberg procedure with arbitrary $q \in [0, 1)$ does not exceed the proportion $1 - \gamma$ of observations that ideally one would like to reject. From a practical point of view, this seems to be a nice "unbiasedness" property of the Benjamini–Hochberg method, one that should be required from procedures for selecting $\mathcal{R}$ in general: at least in the limit, the proportion of false hypotheses among those that pass unnoticed does not exceed the proportion of false hypotheses that would go unnoticed if one simply considered all hypotheses true from the start—if one did not even bother about investigating them—which is just another way of saying that we are better off applying the Benjamini–Hochberg procedure than doing nothing.

For other results on $\Pi_{3,m}$ and a definition of unbiasedness we refer the reader to [16].

**3. Asymptotic results under dependence.** In what follows we assume that

$$(3.1) \qquad F_{m_0}(x) = \frac{1}{m_0} \sum_{k=1}^{m_0} \mathbb{1}_{\{X_k \le x\}} \overset{p}{\to} F(x) := x$$

and

$$(3.2) \qquad G_{m_1}(x) = \frac{1}{m_1} \sum_{k=m_0+1}^{m} \mathbb{1}_{\{X_k \le x\}} \overset{p}{\to} G(x)$$

uniformly in $x \in [0, 1]$, where $G$ is a d.f. concentrated on $[0, 1]$. These are weak versions of the Glivenko–Cantelli theorem; a result at the end of this section gives some sufficient conditions for them to hold.

The following theorem extends Theorem 1 of [9], and in part also Theorem 5 of [19] in the case of the Benjamini–Hochberg method—as opposed to the case of Storey's [18] variant of it (see the Remark to Theorem 4.1 for a parallel result in the case of what we call the Benjamini–Hochberg–Storey method).

THEOREM 3.1. *Under conditions* (3.1) *and* (3.2) *we have, for* $k \in \mathbb{N}$,

$$\underline{\psi}_q^* \left( \frac{q(1-\gamma)}{(1-q)} \right)^k \le \liminf_{m \to \infty} E\left[ \left( \frac{R_m}{m} \right)^k \right]$$

$$\le \limsup_{m \to \infty} E\left[ \left( \frac{R_m}{m} \right)^k \right]$$

$$\le \overline{\psi}_q^* \left( \frac{q(1-\gamma)}{(1-q)} \right)^k,$$

*where, for* $y \ge 0$,

$$\underline{\psi}_q^*(y) = \min\{ x \in [0, 1] : \psi_q(x) \le 1/y \},$$

$$\overline{\psi}_q^*(y) = \inf\{ x \in [0, 1] : \psi_q(x) < 1/y \}$$

*and*

$$\psi_q(x) = \sup_{qx \leq t \leq q} \frac{G(t) - t}{t}, \qquad x \in [0, 1].$$

*In particular,*

(3.3) $$\frac{R_m}{m} \xrightarrow{P} \rho \equiv \rho(q, \gamma) = \psi_q^* \left( \frac{q(1-\gamma)}{(1-q)} \right)$$

*whenever* $\psi_q^*(\frac{q(1-\gamma)}{(1-q)}) := \underline{\psi}_q^*(\frac{q(1-\gamma)}{(1-q)}) = \overline{\psi}_q^*(\frac{q(1-\gamma)}{(1-q)})$, *which will be the case if and only if* $\psi_q$ *does not assume the value* $\frac{(1-q)}{q(1-\gamma)}$ *over an interval.*

PROOF. By (1.2) we have

$$\left\{ \frac{R_m}{m} \geq x \right\} = \left\{ \max_{t = q\lceil mx \rceil/m, \ldots, q(m-1)/m, q} \frac{H_m(t) - t}{t} \geq \frac{1}{q} - 1 \right\}$$

for each $x \in ((r-1)/m, r/m]$, so with $\psi_q^{(m)}(x) = \max_{t = q\lceil mx \rceil/m, \ldots, q} \frac{H_m(t) - t}{t}$,

$$E\left[ \left( \frac{R_m}{m} \right)^k \right] = \int_0^1 k x^{k-1} P\left( \psi_q^{(m)}(x) \geq \frac{1}{q} - 1 \right) dx.$$

Since for each $x > 0$

$$\max_{t = q\lceil mx \rceil/m, \ldots, q(m-1)/m, q} \frac{F_{m_0}(t) - t}{t} \leq \max_{qx \leq t \leq q} \frac{F_{m_0}(t) - t}{t}$$

$$\leq \frac{1}{qx} \max_{qx \leq t \leq q} |F_{m_0}(t) - t| \xrightarrow{P} 0$$

and, similarly, $\max_{t = q\lceil mx \rceil/m, \ldots, q} \frac{G_{m_1}(t) - G(t)}{t} \xrightarrow{P} 0$, we have

$$\psi_q^{(m)}(x) = \gamma \max_{t = q\lceil mx \rceil/m, \ldots, q} \frac{F_{m_0}(t) - t}{t} + (1-\gamma) \max_{t = q\lceil mx \rceil/m, \ldots, q} \frac{G_{m_1}(t) - G(t)}{t}$$

$$+ (1-\gamma) \max_{t = q\lceil mx \rceil/m, \ldots, q} \frac{G(t) - t}{t}$$

$$\xrightarrow{P} (1-\gamma) \max_{qx \leq t \leq q} \frac{G(t) - t}{t}.$$

Thus,

$$\mathbb{1}_{((1-q)/(q(1-\gamma)), \infty)}(\psi_q(x)) \leq \liminf_{m \to \infty} P\left( \psi_q^{(m)}(x) > \frac{1}{q} - 1 \right)$$

$$\leq \limsup_{m \to \infty} P\left( \psi_q^{(m)}(x) \geq \frac{1}{q} - 1 \right)$$

$$\leq \mathbb{1}_{[(1-q)/(q(1-\gamma)), \infty)}(\psi_q(x))$$

for almost all $x$, whence

$$\int_0^1 kx^{k-1} \mathbb{1}_{((1-q)/(q(1-\gamma)),\infty)}(\psi_q(x))\,dx$$

$$\leq \liminf_{m\to\infty} E\left[\left(\frac{R_m}{m}\right)^k\right]$$

$$\leq \limsup_{m\to\infty} E\left[\left(\frac{R_m}{m}\right)^k\right]$$

$$\leq \int_0^1 kx^{k-1} \mathbb{1}_{[(1-q)/(q(1-\gamma)),\infty)}(\psi_q(x))\,dx.$$

Finally, from the definition of $\underline{\psi}_q^*$ and the fact that $\psi_q$ is a nonincreasing right-continuous function, we see that

$$\int_{\{x\in[0,1]\,:\,\psi_q(x)>(1-q)/(q(1-\gamma))\}} kx^{k-1}\,dx = \int_0^{\underline{\psi}_q^*(q(1-\gamma)/(1-q))} kx^{k-1}\,dx$$

$$= \underline{\psi}_q^*\left(\frac{q(1-\gamma)}{(1-q)}\right)^k,$$

the analogous identity for $\overline{\psi}_q^*$ following similarly.  $\square$

REMARKS. (i) Storey, Taylor and Siegmund [19] were the first to realize that conditions like (3.1) and (3.2) are sufficient to derive asymptotic results about Storey's [18] variant of the Benjamini–Hochberg method. Storey, Taylor and Siegmund [19] actually assume only $F(x) \leq x$ in (3.1); assuming $F(x) = x$, however, allows us to obtain simple and useful asymptotic expressions and bounds for $\Pi_{1,m}$, $\Pi_{2,m}$ and $\Pi_{3,m}$ (see the corollaries to the theorem below and Theorem 4.1 later on) without sacrificing much in the domain of practical applicability of the method. Storey, Taylor and Siegmund [19] also assume almost sure convergence in (3.1) and (3.2); our results could as easily be formulated in terms of almost sure convergence, but we find that convergence in probability is more natural in this context—it seems easier to meet and is still very relevant in applications.

(ii) As pointed out by Genovese and Wasserman [9], (3.3) says that asymptotically the Benjamini–Hochberg procedure rejects the observations (or hypotheses whose $p$-values fall) below $q\rho$. Thus, compared with the method of "uncritical multiple testing" in which all hypotheses whose $p$-values fall below a critical value $q$ are rejected, the Benjamini–Hochberg method always rejects a smaller proportion $q\rho(q,\gamma)$ of hypotheses; on the other hand, because $q\rho(q,\gamma) > q/m$ for large $m$, it *typically* rejects many more hypotheses than the corresponding Bonferroni procedure which, for finite $m$, consists of rejecting all observations below $q/m$.

(iii) Suppose (3.3) holds. Then $\rho(q,\gamma) > 0 \Leftrightarrow \max_{qx \le t \le q} \frac{G(t)-t}{t} \ge \frac{(1-q)}{q(1-\gamma)}$ for some $x > 0$, and it can be seen that

$$(3.4) \qquad \rho(q,\gamma) = q^{-1} \sup\left\{x \in [0,1] : \frac{G(x)-x}{x} > \frac{(1-q)}{q(1-\gamma)}\right\},$$

or $q\rho(q,\gamma) = \sup\{x \in [0,1] : \frac{x}{H(x)} < q\}$, in agreement with Theorem 5 of Storey, Taylor and Siegmund [19]. Furthermore, it can be verified from (3.4) that $\rho(q,\gamma)$ is left-continuous in $q$ for fixed $\gamma$, and, using the condition expressed right after (3.3), that it is right-continuous at $q$ if $\rho(q,\gamma) > 0$. Thus, $q \to \rho(q,\gamma)$ is continuous on $(q',q'')$ if $\rho(q,\gamma) > 0 \;\forall q \in (q',q'')$, in which case $R_m/m \overset{P}{\to} \rho(q,\gamma)$ uniformly on $[q',q'']$. ($R_m/m$ is a nondecreasing right-continuous function of $q$.)

EXAMPLES. (i) Suppose $G$ is degenerate at $x_0 \in [0,1)$. Then $\psi_q(x) = -1$ if $q < x_0$, $\psi_q(x) = 1/x_0 - 1$ if $qx \le x_0 < q$, and $\psi_q(x) = 1/x - 1$ if $qx \ge x_0$.

If $x_0 > q(1-\gamma)/(1-q\gamma)$, that is, if $1/x_0 - 1 < (1-q)/[q(1-\gamma)]$, then $\overline{\psi_q^*}\left(\frac{q(1-\gamma)}{(1-q)}\right) = 0$, and hence $\rho = 0$.

If $x_0 < q(1-\gamma)/(1-q\gamma)$, then the equation $\psi_q(x) = (1-q)/[q(1-\gamma)]$ has a unique solution given by $x = (1-\gamma)/(1-q\gamma)$, so (3.3) holds and

$$(3.5) \qquad \rho(q,\gamma) = \psi_q^*\left(\frac{q(1-\gamma)}{(1-q)}\right) = \frac{(1-\gamma)}{(1-q\gamma)}.$$

Thus, $\rho(q,\gamma) > 0$ if $x_0 < (1-\gamma)/(1-q\gamma)$, that is, $\rho > 0$ if $x_0$ is not "too large" given the choice of $q$, in which case $\rho$ is actually independent of $x_0$, implying that asymptotically the proportion of rejections and the efficacy of the procedure depend only on $\gamma$ and on the choice of $q$ and not on the exact position of $x_0$. In fact, it can be checked by substitution of (3.5) into the expressions of the limits obtained below in (3.8) that $\Pi_{2,m}$ and $1-\Pi_{3,m}$ both converge in probability to 1 when $x_0 < (1-\gamma)/(1-q\gamma)$.

Since $q$ can always be chosen so that $x_0 < (1-\gamma)/(1-q\gamma)$, we see that in this case the Benjamini–Hochberg procedure can always be made to work in an asymptotically optimal way—in such a way that practically 100% of the observations from $G$ will be spotted and $\Pi_{1,m}$ is kept at $q\gamma$. In order to make use of this optimality in practice, one needs to choose $q$ appropriately, but this is easy if $\gamma$ is not too large, because the histogram will then have the shape of a scaled down uniform density with a conspicuous peak at $x_0$ (which is why the problem is easy to solve even without using the Benjamini–Hochberg method).

In the borderline case where $x_0 = q(1-\gamma)/(1-q\gamma)$, the theorem only tells us that $R_m/m$ is asymptotically somewhere between 0 and the right-hand side of (3.5), because $\psi_q(x) = (1-q)/[q(1-\gamma)] \equiv x_0/q$ holds for all

$x \in (0, x_0/q)$. In fact, *if $X_1, \ldots, X_{m_0}$ are independent standard uniform* r.v.'s, we have

$$(3.6) \qquad \frac{R_m}{m} \xrightarrow{P} \frac{(1-\gamma)}{2(1-q\gamma)}.$$

To see this, note that, after being sorted in ascending order, the sample consists of a proportion $H_m(x_0-)$ of ordered uniforms below $x_0$, followed by $m - m_0$ copies of $x_0$, which are in turn followed by the remaining $m(1 - H_m(x_0))$ ordered uniforms, so that the proportion of correctly rejected observations is always given by $(R_m - S_m)/m = \max\{i : m_0 H_m(x_0-) < i \le m - m(1 - H_m(x_0)), mx_0/q \le i\}/m - H_m(x_0-)$. This is $\neq 0$ and equals $1 - \gamma$ if and only if $mx_0/q \le m - m(1 - H_m(x_0))$, or $F_{m_0}(x_0) - x_0 \ge 0$, which by our assumption happens with probability tending to $1/2$. Thus, $\frac{R_m - S_m}{m} \xrightarrow{P} \frac{1-\gamma}{2}$, and therefore (3.6) holds by the fact that $S_m/R_m \to q\gamma$, which follows by Proposition 2.2 (note that $R_m \xrightarrow{P} \infty$ necessarily).

Finally, we observe that in this borderline case $\Pi_{2,m}$ and $1 - \Pi_{3,m}$ converge in probability to $1/2$ and $1 - (1-\gamma)(1-q\gamma)/[(1-q\gamma) + \gamma(1-q)]$, respectively, a calculation suggesting that $\Pi_{2,m}$ is a more practically meaningful measure of efficacy than $1 - \Pi_{3,m}$.

(ii) Assume that $G$ is concave and

$$G'_+(0) = \lim_{x \downarrow 0} \frac{G(x)}{x} > \beta \qquad \text{where } \beta = \frac{(1-q\gamma)}{q(1-\gamma)}.$$

Since then $G(0) = 0$ and $\beta \ge 1$, there exists a unique $t^* > 0$ such that $G(t^*) = \beta t^*$; moreover, $t^* \le q$ [because $1 \ge G(t^*) = \beta t^* = \frac{t^*}{q}\frac{1-q\gamma}{1-\gamma}$ and $\frac{1-q\gamma}{1-\gamma} \ge 1 \Rightarrow t^* \le q$], and it becomes evident on geometric grounds that

$$\max_{t^* \le t \le q} \frac{G(t)}{t} - 1 = \frac{G(t^*)}{t^*} - 1 = \beta - 1 < \max_{qx \le t \le q} \frac{G(t)}{t} - 1 \qquad \forall x \in (0, t^*/q);$$

thus,

$$\rho(q, \gamma) = \psi_q^* \left( \frac{q(1-\gamma)}{(1-q)} \right) = \psi_q^* \left( \frac{1}{\beta - 1} \right) = \frac{t^*}{q}.$$

Alternatively, by (3.4), $q\rho(q, \gamma)$ is the smallest positive root of $G(t) = \beta t$, that is, $q\rho(q, \gamma) = t^*$. This was first proved by Genovese and Wasserman [9].

(iii) For an example where $G$ is not necessarily concave take $G(x) = px^\alpha + (1-p)x^\beta$, $0 \le x \le 1$, with $\alpha \in (0,1)$, $\beta > 1$, $0 < p < 1$. Then $(G(t) - t)/t = pt^{\alpha-1} + (1-p)t^{\beta-1} - 1$, and from (3.4) we see that $\rho > 0$ always exists and is uniquely determined by $p(q\rho)^{\alpha-1} + (1-p)(q\rho)^{\beta-1} - 1 = \frac{(1-q)}{q(1-\gamma)}$, provided $q > 0$.

Using Theorem 3.1, we can show that the conclusion of the Benjamini–Hochberg theorem holds very generally in an asymptotic sense:

COROLLARY 3.2.    *Under the conditions of Theorem* 3.1,

$$(3.7) \qquad \frac{R_m}{m} \xrightarrow{P} \rho > 0 \quad \Longrightarrow \quad \Pi_{1,m} \xrightarrow{P} q\gamma.$$

PROOF.   Since

$$\Pi_{1,m} = \frac{S_m}{R_m \vee 1} = \frac{\sum_{i=1}^{m_0} \mathbb{1}_{\{X_i \leq qR_m/m\}}}{R_m \vee 1} = \gamma \frac{(1/m_0)\sum_{i=1}^{m_0} \mathbb{1}_{\{X_i \leq qR_m/m\}}}{(R_m \vee 1)/m},$$

we have for arbitrary $\varepsilon \in (0, \rho)$, $\eta \in (0, 1)$,

$$\gamma \frac{F_{m_0}(q(\rho - \varepsilon))}{\rho + \varepsilon} \leq \Pi_{1,m} \leq \gamma \frac{F_{m_0}(q(\rho + \varepsilon))}{\rho - \varepsilon},$$

with probability at least $1 - \eta$, which by (3.1) proves (3.7).   □

The following statements are all direct consequences of the preceding results.

COROLLARY 3.3.    *Under the conditions of Theorem* 3.1, $\frac{R_m}{m} \xrightarrow{P} \rho(q, \gamma) > 0$ *implies*

$$(3.8) \qquad \begin{aligned} &\frac{S_m}{m} \xrightarrow{P} \rho(q, \gamma) q\gamma, \\ &\Pi_{2,m} \xrightarrow{P} \rho(q, \gamma)\frac{(1 - q\gamma)}{(1 - \gamma)} \quad \text{and} \quad 1 - \Pi_{3,m} \xrightarrow{P} \gamma \frac{(1 - q\rho(q, \gamma))}{(1 - \rho(q, \gamma))}. \end{aligned}$$

Because $R_m/m$, $S_m/m$, $\Pi_{1,m}$, $\Pi_{2,m}$ and $\Pi_{3,m}$ are proportions, all the above statements about convergence in probability to a constant are equivalent to statements about convergence in the mean (of any order), as well as to statements about convergence of their moments. One consequence of this fact is that, *under the conditions of Theorem* 3.1,

$$\frac{R_m}{m} \xrightarrow{P} \rho > 0 \quad \Longrightarrow \quad \lim_{m \to \infty} \frac{E[S_m{}^k]}{E[R_m{}^k]} = \lim_{m \to \infty} E\left[\left(\frac{S_m}{R_m}\right)^k\right] = (q\gamma)^k,$$

which implies that, asymptotically, the Benjamini–Hochberg method also controls the "marginal false discovery rate" $E(S_m)/E(R_m \vee 1)$ [briefly mentioned in Remark (ii) to Proposition 2.2].

We shall finish this section by giving an example of a rather general situation in which statements like (3.1) and (3.2) hold true uniformly in $x$; a similar result (with a stronger conclusion) for stationary ergodic sequences has been given by Tucker [21], for example. Let $\xi_1, \xi_2, \ldots$ be a sequence of r.v.'s on $[0, 1]$ with d.f.'s $G^{(1)}, G^{(2)}, \ldots$. Since for each $x$ $G_n(x) := n^{-1}\sum_{i=1}^{n} \mathbb{1}_{\{\xi_i \leq x\}} \xrightarrow{P}$ $G(x)$ if and only if $EG_n(x) \to G(x)$ and $E(G_n(x)^2) \to G(x)^2$, we see that $G_n(x) \xrightarrow{P} G(x)$ is equivalent to

$$\lim_{n\to\infty} \frac{1}{n} \sum_{i=1}^{n} G^{(i)}(x) = G(x) \quad \text{and} \quad \lim_{n\to\infty} \frac{1}{n^2} \sum_{i\neq j}^{n} P(\xi_i \leq x, \xi_j \leq x) = G(x)^2.$$

The following sufficient condition combines this observation with a condition that is much weaker than strong mixing.

PROPOSITION 3.4. *Assume that, for each $x$,*

$$G(x) := \lim_{n\to\infty} \frac{1}{n} \sum_{i=1}^{n} G^{(i)}(x) \quad and \quad G(x-) := \lim_{n\to\infty} \frac{1}{n} \sum_{i=1}^{n} G^{(i)}(x-)$$

*exist, and there are subsequences $\{k_n\}$ and $\{\alpha_{k_n}\}$ such that $k_n \to \infty$, $k_n/n \to 0$ and $\alpha_{k_n} \to 0$ as $n \to \infty$, and*

$$\sup_{|i-j|\geq k_n} \max\{|P(\xi_i \leq x, \xi_j \leq x) - P(\xi_i \leq x)P(\xi_j \leq x)|,$$

$$|P(\xi_i < x, \xi_j < x) - P(\xi_i < x)P(\xi_j < x)|\} \leq \alpha_{k_n}.$$

*Then $G_n \xrightarrow{P} G$ uniformly.*

PROOF. That $G_n(x) \xrightarrow{P} G(x)$ for fixed $x$ follows from the fact that $\lim_{n\to\infty} \frac{1}{n^2} \sum_{i\neq j}^{n} P(\xi_i \leq x)P(\xi_j \leq x) = \lim_{n\to\infty} (\frac{1}{n} \sum_{i=1}^{n} P(\xi_i \leq x))^2 = G(x)^2$ and from the inequalities

$$\left| \frac{1}{n^2} \sum_{i\neq j}^{n} P(\xi_i \leq x, \xi_j \leq x) - \frac{1}{n^2} \sum_{i\neq j}^{n} P(\xi_i \leq x)P(\xi_j \leq x) \right|$$

$$\leq \left( \frac{k_n}{n} \right)^2 + \frac{1}{n^2} \sum_{|i-j|\geq k_n}^{n} |P(\xi_i \leq x, \xi_j \leq x) - P(\xi_i \leq x)P(\xi_j \leq x)|$$

$$\leq \left( \frac{k_n}{n} \right)^2 + \alpha_{k_n}$$

(the right-hand side of which goes to zero as $n \to \infty$ by assumption). The analogous statement with $< x$ in place of $\leq x$ and $x-$ in place of $x$ follows in the same way. Finally, that these pointwise results imply uniform convergence is a classical result.  $\square$

**4. A modification of the method.** It has been observed by several authors that the Benjamini–Hochberg method tends to be conservative unless $\gamma$ is relatively close to 1. For if the value of $\gamma$ cannot be guessed at, the only way one can guarantee that $E(\Pi_{1,m}) \leq \delta$ for a given $\delta > 0$ is to apply the method with $q = \delta$. But if $\gamma$ is actually smaller, say equal to 1/2, such a choice yields the overcautious bound $E(\Pi_{1,m}) \leq \delta/2$ and the concomitant decrease in $\Pi_{2,m}$, which is an increasing function of $q$. Although in some practical situations this is hardly a problem because one has a reasonably good idea about the value of $\gamma$, from a general point of view it is still a shortcoming one would like to eliminate.

These considerations have led Benjamini and Hochberg [3], Storey [18] and Storey, Taylor and Siegmund [19], among others, to propose and study variants of the Benjamini–Hochberg method which incorporate estimates of $\gamma$. Our objective here will be to introduce another variant—very similar to Storey's—and to study some of its asymptotic properties. Questions related to the practical application of the method [e.g., the problem of choosing $x$ in (4.1) below] will be considered elsewhere. Our assumptions and notation will be those of Section 3.

The closer $x$ gets to $G^{-1}(1)$, the tighter the inequality $H(x) = \gamma x + (1 - \gamma)G(x) \leq \gamma x + (1 - \gamma)$, or $\gamma \leq \frac{1 - H(x)}{1 - x}$, becomes, which suggests taking

$$(4.1) \qquad \gamma_m(x) = \min_{0 \leq t \leq x} \frac{1 - H_m(t)}{1 - t},$$

where $x \in (0, 1)$ is to be chosen, as an estimator of $\gamma$ [note that, for fixed $x \in (0, 1)$, $\gamma_m(x) > 0$ with probability tending to 1]. (Storey's [18] estimator is defined by $(1 - H_m(x))/(1 - x)$ for a given $x$.) Because of the convergence of $H_m$ to $H$, this $\gamma_m(x)$ will typically be an overestimate of $\gamma$ in the sense that, given $\varepsilon > 0$,

$$(4.2) \qquad \gamma_m(x) = \min_{0 \leq t \leq x} \frac{1 - H_m(t)}{1 - t} > \min_{0 \leq t \leq x} \frac{1 - H(t)}{1 - t} - \varepsilon \geq \gamma - \varepsilon,$$

with high probability if $m$ is large enough. On the other hand, if we put

$$\kappa(x) = \min_{0 \leq t \leq x} \frac{1 - G(t)}{1 - t}, \qquad x \in (0, 1),$$

we see that $\gamma_m(x)$ will typically not exceed $\gamma$ by more than $(1 - \gamma)\kappa(x)$:

$$(4.3) \qquad \begin{aligned} \gamma_m(x) &= \gamma \min_{0 \leq t \leq x} \frac{1 - F_{m_0}(t)}{1 - t} + (1 - \gamma) \min_{0 \leq t \leq x} \frac{1 - G_{m_1}(t)}{1 - t} \\ &< \varepsilon + \gamma + (1 - \gamma)\kappa(x), \end{aligned}$$

with high probability for arbitrary $\varepsilon > 0$ if $m$ is large enough.

For want of a better name, and because we are essentially using the ideas of Benjamini and Hochberg [2] and Storey [18], we shall refer to the procedure that consists of rejecting all observations smaller than or equal to $X_{R_m(q_m(x,\delta)):m}$, where $R_m(q_m(x,\delta)) = \max\{i : X_{i:m} \leq q_m(x,\delta)\frac{i}{m}\}$, $q_m(x,\delta) = \frac{\delta}{\gamma_m(x)}$ and $\gamma_m(x)$ is defined by (4.1), as the *Benjamini–Hochberg–Storey method*.

The variable $R_m$ of (1.1) will now be denoted by $R_m(q)$ in order to indicate its dependence on $q$ in the Benjamini–Hochberg method, and similarly for the other variables; for instance, we shall write $\Pi_{1,m}(q)$ for $\Pi_{1,m}$, and $\Pi_{1,m}(q_m(x,\delta))$ for the proportion of incorrect rejections incurred by applying the Benjamini–Hochberg–Storey method.

The following result shows that, with the modified method, one is able, in an asymptotic sense, to keep the false discovery rate under control and at the same time achieve greater average power than that provided by the Benjamini–Hochberg procedure.

THEOREM 4.1. *Let* $\gamma \in (0,1)$ *and suppose* $\delta > 0$, $x \in (0,1)$, $q'(x)$ *and* $q''(x)$ *can be chosen so that*

$$q'(x) < \frac{\delta}{\gamma + (1-\gamma)\kappa(x)} \leq \frac{\delta}{\gamma} < q''(x)$$

*and*

$$\frac{R_m(q)}{m} \xrightarrow{P} \rho(q,\gamma) > 0 \qquad \forall q \in [q'(x), q''(x)].$$

*Then*

(4.4)
$$\delta \frac{\gamma}{\gamma + (1-\gamma)\kappa(x)} \leq \liminf_{m\to\infty} E[\Pi_{1,m}(q_m(x,\delta))]$$
$$\leq \limsup_{m\to\infty} E[\Pi_{1,m}(q_m(x,\delta))] \leq \delta$$

*and*

(4.5)
$$\rho\left(\frac{\delta}{\gamma + (1-\gamma)\kappa(x)}, \gamma\right) \frac{1 - \delta\gamma/(\gamma + (1-\gamma)\kappa(x))}{1-\gamma}$$
$$\leq \liminf_{m\to\infty} E[\Pi_{2,m}(q_m(x,\delta))]$$
$$\leq \limsup_{m\to\infty} E[\Pi_{2,m}(q_m(x,\delta))] \leq \rho\left(\frac{\delta}{\gamma}, \gamma\right) \frac{(1-\delta)}{(1-\gamma)}.$$

PROOF. We know from Corollary 3.3 that we have

$$\frac{R_m(q)}{m} \xrightarrow{P} \rho(q,\gamma) \qquad \text{as well as} \qquad \frac{S_m(q)}{m} \xrightarrow{P} \rho(q,\gamma)q\gamma$$

$\forall q \in [q'(x), q''(x)]$; moreover, by Remark (iii) following Theorem 3.1, the convergence here is uniform on $[q'(x), q''(x)]$. It can be shown (and it is certainly known) that if $f_n \to f$ and $g_n \to g$ uniformly, $\sup_t |f(t)| < \infty$ and $\inf_t |g(t)| > 0$, then $\sup_t |f_n(t)/g_n(t) - f(t)/g(t)| \to 0$. Thus,

$$(4.6) \qquad \sup_{q'(x) \le q \le q''(x)} \left| \frac{S_m(q)}{R_m(q) \vee 1} - q\gamma \right| = \sup_{q'(x) \le q \le q''(x)} |\Pi_{1,m}(q) - q\gamma| \xrightarrow{P} 0.$$

Now fix $\varepsilon \in (0, \gamma)$, $\eta \in (0, 1)$ and $m'$ so large that

$$(4.7) \qquad \begin{aligned} q'(x) &\le \frac{\delta}{\gamma + (1-\gamma)\kappa(x) + \varepsilon} \\ &\le q_m(x, \delta) \equiv \frac{\delta}{\gamma_m(x)} \le \frac{\delta}{\gamma - \varepsilon} \le q''(x), \end{aligned}$$

with probability at least $1 - \eta$ if $m \ge m'$, which is possible by (4.2), (4.3) and our assumptions about $q'(x)$ and $q''(x)$. Then for $m \ge m'$

$$\Pi_{1,m}(q_m(x, \delta)) \le \sup_{q'(x) \le q \le q''(x)} |\Pi_{1,m}(q) - q\gamma| + \delta \frac{\gamma}{\gamma - \varepsilon}$$

holds with probability at least $1 - \eta$. Since $\varepsilon$ is arbitrarily small, this, combined with (4.6), proves the inequality on the right-hand side in (4.4) as well as its version in probability. The other inequality follows similarly.

To prove (4.5), we use the inequalities

$$\frac{R_m(\delta/(\gamma + (1-\gamma)\kappa(x) + \varepsilon)) - S_m(\delta/(\gamma + (1-\gamma)\kappa(x) + \varepsilon))}{m - m_0}$$

$$\le \frac{R_m(q_m(x, \delta)) - S_m(q_m(x, \delta))}{m - m_0}$$

$$\le \frac{R_m(\delta/(\gamma - \varepsilon)) - S_m(\delta/(\gamma - \varepsilon))}{m - m_0},$$

which hold whenever (4.7) is valid because $R_m(q) - S_m(q)$ is nondecreasing in $q$, and the continuity of $q \to \rho(q, \gamma)$ on $[q'(x), q''(x)]$. □

REMARK. Under the assumptions of the theorem, we have $q_m(x, \delta) \xrightarrow{P} q(x, \delta) := \frac{\delta}{\gamma + (1-\gamma)\kappa(x)}$ and $\frac{R_m(q_m(x, \delta))}{m} \xrightarrow{P} \rho(q(x, \delta), \gamma)$; thus, asymptotically, the Benjamini–Hochberg–Storey method consists of rejecting all observations below $q(x, \delta)\rho(q(x, \delta), \gamma)$.

EXAMPLES. (i) If $G(x) = x^\alpha$, $x \in [0, 1]$, $\alpha \in (0, 1)$, then $\kappa(x) = (1 - x^\alpha)/(1 - x)$ because $t \to (1 - t^\alpha)/(1 - t)$ is decreasing. [In fact, if $G$ has a nonincreasing density function $g$, then $1 - G(t) = \int_t^1 g(s)\,ds \le (1-t)g(t)$, or $-g(t)(1-t) + (1 - G(t)) < 0$, which implies that the derivative of $t \to (1 - G(t))/(1-t)$ is negative.] In this case [see Example (ii) following Theorem 3.1], it can be seen that $\rho(q, \gamma) = (q(1-\gamma)/(1 - q\gamma))^{1/(1-\alpha)}/q$, which is always positive for $q > 0$, and so we have explicit expressions for the bounds in Theorem 4.1 that are valid for all $x \in (0, 1)$. Here we shall consider $\alpha = 0.1$ in two cases: (a) $\gamma = 0.5$, (b) $\gamma = 0.9$. The density $h$ of $H$ in case (a) is roughly in agreement with the histogram shown in Figure 5.8 of [12]; that of case (b) is much closer to the standard uniform density; they are both compared with the latter in Figure 1.

The asymptotic average power and false discovery rate of the Benjamini–Hochberg procedure are shown in Figure 2 as functions of $q$. In case (a), the choice of $q = 0.2$ yields an asymptotic false discovery rate of 0.1 and an asymptotic average power of 0.784; in case (b), an asymptotic false discovery rate of 0.1 is guaranteed by taking $q = 0.111$, which yields an asymptotic average power of 0.614.

Figure 3 illustrates the adherence of the bounds in (4.5) as a function of $x$ when $\delta$ (the upper bound of the false discovery rate) is fixed at 0.1; as just seen, in the ideal situation where $\gamma$ is known, the power obtained by controlling the false discovery rate at this level would be about 0.784 and 0.614 in the cases $\gamma = 0.5$ and $\gamma = 0.9$, respectively. In each case, the asymptotic average power of the Benjamini–Hochberg–Storey procedure with $q_m(x, \delta) = 0.1/\gamma_m(x)$ lies between the two curves of Figure 3 and is rather close to the maximum average power—achieved by setting $q = \delta/\gamma$ in the Benjamini–Hochberg procedure—even for small values of $x$. However, since $\kappa(x) \to \alpha$ as $x \uparrow 1$, the lower bound for asymptotic average power is always strictly below $\rho(\delta/(\gamma + (1-\gamma)\alpha), \gamma)(1 - \delta\frac{\gamma}{\gamma + (1-\gamma)\alpha})/(1 - \gamma)$, which in turn is always strictly below the asymptotic average power of the Benjamini–Hochberg procedure with $q = \delta/\gamma$.

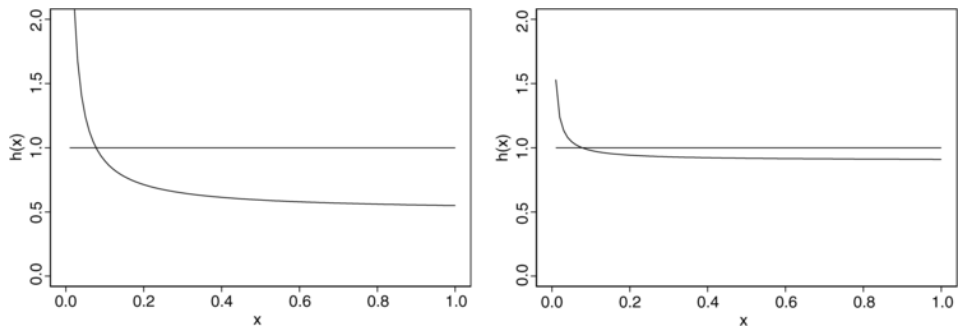

FIG. 1. *Densities of the standard uniform distribution and of the d.f. $H$: left panel:* $\alpha = 0.1$, $\gamma = 0.5$, *right panel:* $\alpha = 0.1$, $\gamma = 0.9$.

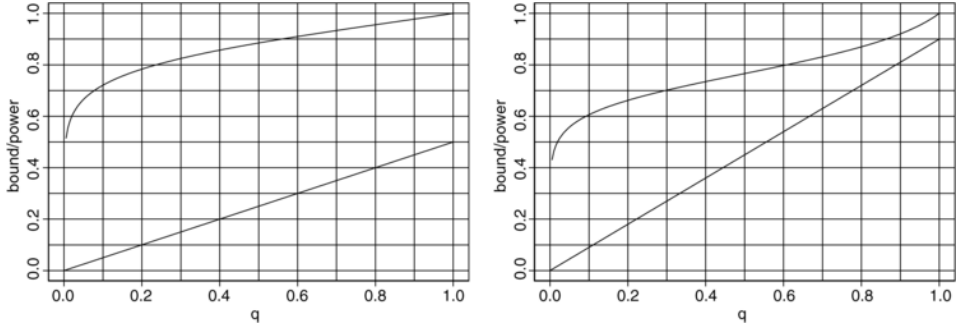

Fig. 2. *Asymptotic average power and false discovery rate of the Benjamini–Hochberg procedure as functions of $q$:* left panel: $\alpha = 0.1$, $\gamma = 0.5$, right panel: $\alpha = 0.1$, $\gamma = 0.9$.

The left-hand side of (4.4) approaches $\delta = 0.1$ in a very similar way.

(ii) Suppose $G(x) = x^{\alpha}\mathbb{1}_{[0,x_0)}(x) + \mathbb{1}_{[x_0,\infty)}(x)$ for $x_0, \alpha \in (0,1)$. Then $\psi_q(x) = (qx)^{\alpha-1} - 1$ if $0 \leq x < x_0/q$ and $\psi_q(x) = 0$ if $x \geq x_0/q$, so that $\rho(q, \gamma)$ is still positive and has the same expression as in (i) as long as $(q(1-\gamma)/(1-q\gamma))^{1/(1-\alpha)} \leq x_0$, which can always be arranged by choosing a small enough $q$. Since $\kappa(x) = (1-x^{\alpha})/(1-x)$ for $x \in [0, x_0)$ and $\kappa(x) = 0$ for $x \in [x_0, 1)$, the lower bounds on the average power of the Benjamini–Hochberg–Storey procedure as a function of $x$ coincide in this case with those shown in Figure 3 over the interval $[0, x_0)$, but attain their maximum values over $[x_0, 1)$; analogously, the lower bounds on the false discovery rate attain the value of $\delta$ if $x \in [x_0, 1)$.

In this case, therefore, using $q_m(x, \delta) = \delta/\gamma_m(x)$ with $x \in [x_0, 1)$ in place of $q$ in the Benjamini–Hochberg procedure and choosing $\delta$ according to the conditions of Theorem 4.1 is asymptotically equivalent to taking $q = \delta/\gamma$ and thus corresponds to the ideal situation in which $\gamma$ is known, the required upper bound for the false discovery rate is $\delta$, and the power is maximum.

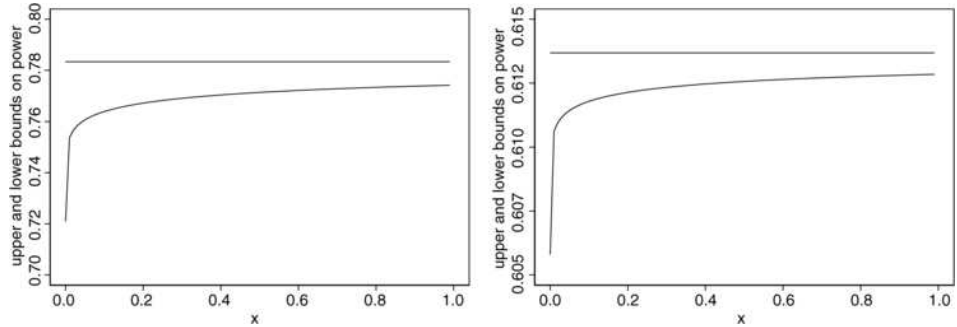

Fig. 3. *Upper and lower bounds on the asymptotic average power of the Benjamini–Hochberg–Storey procedure as functions of $x$ as given in* (4.5)*:* left panel: $\alpha = 0.1$, $\gamma = 0.5$, right panel: $\alpha = 0.1$, $\gamma = 0.9$.

## APPENDIX: PROOF OF PROPOSITION 2.2

We first show that $R_m \xrightarrow{P} \infty \Rightarrow R_m^{(j)}(X^{(j)}) \xrightarrow{P} \infty \ \forall j$. Observe that $H_m^*(t) := m^{-1} \sum_{i=2}^m \mathbb{1}_{\{X_i \le t + q/m\}} \ge \widetilde{H}_m(t) := m^{-1} \sum_{i=2}^m \mathbb{1}_{\{X_i \le t\}}$ for all $t$, and that, by definition of $R_m^{(1)}(X^{(1)})$ and (1.2), we have

$$\left\{ \frac{R_m^{(1)}(X^{(1)})}{m-1} \ge x \right\} = \left\{ \max_{t=qr/m,\ldots,q(m-1)/m} \frac{H_m^*(t) - t}{t} \ge \frac{1}{q} - 1 \right\}$$

for $x \in ((r-1)/(m-1), r/(m-1)]$. Since

$$\max_{t=qr/m,\ldots,q(m-1)/m,q} \frac{\widetilde{H}_m(t) - t}{t}$$

$$= \max\left\{ \max_{t=qr/m,\ldots,q(m-1)/m} \frac{\widetilde{H}_m(t) - t}{t}, \frac{\widetilde{H}_m(q) - q}{q} \right\}$$

$$\le \max\left\{ \max_{t=qr/m,\ldots,q(m-1)/m} \frac{H_m^*(t) - t}{t}, \right.$$

$$\left. \frac{H_m^*(q(m-1)/m) - q(m-1)/m}{q(m-1)/m} \frac{(m-1)}{m} - \frac{1}{m} \right\},$$

we have

$$\frac{1}{q} - 1 \le \max_{t=qr/m,\ldots,q(m-1)/m,q} \frac{\widetilde{H}_m(t) - t}{t}$$

$$\implies \frac{1}{q} - 1 \le \max_{t=qr/m,\ldots,q(m-1)/m} \frac{H_m^*(t) - t}{t},$$

and because $\sup_t |\widetilde{H}_m(t) - H_m(t)| \to 0$ with probability one (and $q\frac{r}{m} > qx/2$), it follows that $P(R_m/m \ge x(m-1)/m) + \varepsilon \le P(R_m^{(1)}(X^{(1)})/(m-1) \ge x)$ for sufficiently large $m$ and arbitrary $\varepsilon > 0$. This proves that $R_m \xrightarrow{P} \infty \Rightarrow R_m^{(1)}(X^{(1)}) \xrightarrow{P} \infty$; similar reasoning shows that $R_m^{(j)}(X^{(j)}) \xrightarrow{P} \infty \Rightarrow R_m^{(j+1)}(X^{(j+1)}) \xrightarrow{P} \infty$. Thus, $R_m \xrightarrow{P} \infty$ implies $R_m^{(j)}(X^{(j)}) \xrightarrow{P} \infty$ for each $j$, and by the bounded convergence theorem $E[(j + R_m^{(j)}(X^{(j)}))^{j-k}] \to 0$ whenever $1 \le j < k$, so (2.3) follows from (2.1). In the standard uniform case equality holds in (2.3) with "lim" in place of "lim sup," whence $\Pi_{1,m} \xrightarrow{P} q\gamma$.

To prove the converse, we show that $\Pi_{1,m} \xrightarrow{P} q\gamma \Rightarrow R_m^{(1)}(X^{(1)}) \xrightarrow{P} \infty$ and then that $R_m^{(1)}(X^{(1)}) \xrightarrow{P} \infty \Rightarrow R_m \xrightarrow{P} \infty$. Suppose $\Pi_{1,m} \xrightarrow{P} q\gamma$, and assume $\limsup_{m \to \infty} R_m^{(1)}(X^{(1)}) \leq C < \infty$ in probability. Then (2.1) with $k = 2$ and in the standard uniform case implies

$$\liminf_{m \to \infty} E[\Pi_{1,m}^2] = q\gamma \liminf_{m \to \infty} E\left[\frac{1}{1 + R_m^{(1)}(X^{(1)})}\right] + (q\gamma)^2$$

$$\geq \frac{q\gamma}{1 + C} + (q\gamma)^2 > (q\gamma)^2,$$

which contradicts $\Pi_{1,m} \xrightarrow{P} q\gamma$; thus, $R_m^{(1)}(X^{(1)}) \xrightarrow{P} \infty$. When $k = 1$, (2.2) in the standard uniform case reads

$$(A.1) \qquad \frac{E(S_m)}{1 + E(R_m^{(1)}(X^{(1)}))} = q\gamma.$$

If $R_m \xrightarrow{P} \infty$ then $S_m \xrightarrow{P} \infty$, but then $R_m^{(1)}(X^{(1)}) \xrightarrow{P} \infty$ contradicts (A.1) when we let $m \to \infty$; thus we must have $R_m \xrightarrow{P} \infty$ if $R_m^{(1)}(X^{(1)}) \xrightarrow{P} \infty$.

## REFERENCES

[1] Abramovich, F., Benjamini, Y., Donoho, D. and Johnstone, I. M. (2006). Adapting to known sparsity by controlling the false discovery rate. *Ann. Statist.* **34** 584–653.

[2] Benjamini, Y. and Hochberg, Y. (1995). Controlling the false discovery rate: A practical and powerful approach to multiple testing. *J. Roy. Statist. Soc. Ser. B* **57** 289–300. MR1325392

[3] Benjamini, Y. and Hochberg, Y. (2000). On the adaptive control of the false discovery rate in multiple testing with independent statistics. *J. Educational and Behavioral Statistics* **25** 60–83.

[4] Benjamini, Y. and Yekutieli, D. (2001). The control of the false discovery rate in multiple testing under dependence. *Ann. Statist.* **29** 1165–1188. MR1869245

[5] Delongchamp, R., Bowyer, J., Chen, J. and Kodell, R. (2004). Multiple-testing strategy for analyzing cDNA array data on gene expression. *Biometrics* **60** 774–782. MR2089454

[6] Donoho, D. and Jin, J. (2004). Higher criticism for detecting sparse heterogeneous mixtures. *Ann. Statist.* **32** 962–994. MR2065195

[7] Dudoit, S., Schaffer, J. and Boldrick, J. (2003). Multiple hypothesis testing in microarray experiments. *Statist. Sci.* **18** 71–103. MR1997066

[8] Fan, J., Tam, P., Vande Woude, G. and Ren, Y. (2004). Normalization and analysis of cDNA microarrays using within-array replications applied to neuroblastoma cell response to a cytokine. *Proc. Natl. Acad. Sci. USA* **101** 1135–1140.

[9] Genovese, C. and Wasserman, L. (2002). Operating characteristics and extensions of the false discovery rate procedure. *J. R. Stat. Soc. Ser. B Stat. Methodol.* **64** 499–517. MR1924303

[10] Genovese, C. and Wasserman, L. (2004). A stochastic process approach to false discovery control. *Ann. Statist.* **32** 1035–1061. MR2065197

[11] Karlin, S. and Taylor, H. (1981). *A Second Course in Stochastic Processes.* Academic Press, New York. MR0611513

[12] McLachlan, G., Do, K.-A. and Ambroise, C. (2004). *Analyzing Microarray Gene Expression Data.* Wiley, New York.

[13] Nakas, C., Yiannoutsos, C., Bosch, R. and Moyssiadis, C. (2003). Assessment of diagnostic markers by goodness-of-fit tests. *Statistics in Medicine* **22** 2503–2513.

[14] Reiner, A., Yekutieli, D. and Benjamini, Y. (2003). Identifying differentially expressed genes using false discovery rate controlling procedures. *Bioinformatics* **19** 368–375.

[15] Sarkar, S. (2002). Some results on false discovery rate in stepwise multiple testing procedures. *Ann. Statist.* **30** 239–257. MR1892663

[16] Sarkar, S. (2004). FDR-controlling stepwise procedures and their false negative rates. *J. Statist. Plann. Inference* **125** 119–137. MR2086892

[17] Sarkar, S. (2006). False discovery and false non-discovery rates in single-step multiple testing procedures. *Ann. Statist.* **34** 394–415.

[18] Storey, J. (2002). A direct approach to false discovery rates. *J. R. Stat. Soc. Ser. B Stat. Methodol.* **64** 479–498. MR1924302

[19] Storey, J., Taylor, J. and Siegmund, D. (2004). Strong control, conservative point estimation and simultaneous conservative consistency of false discovery rates: A unified approach *J. R. Stat. Soc. Ser. B Stat. Methodol.* **66** 187–205. MR2035766

[20] Tsai, C.-A., Hsueh, H.-M. and Chen, J. (2003). Estimation of false discovery rates in multiple testing: Application to gene microarray data. *Biometrics* **59** 1071–1081. MR2025132

[21] Tucker, H. (1953). A generalization of the Glivenko–Cantelli theorem. *Ann. Math. Statist.* **30** 828–830. MR0107891

[22] Tusher, V., Tibshirani, R. and Chu, G. (2001). Significance analysis of microarrays applied to the ionizing radiation response. *Proc. Natl. Acad. Sci. USA* **98** 5116–5121.

Department of Clinical Epidemiology
and Biostatistics
Academic Medical Centre
University of Amsterdam
P.O. Box 22700
1100 DE Amsterdam
The Netherlands
E-mail: J.A.Ferreira@amc.uva.nl
          ahzwinderman@amc.uva.nl